\documentclass[11pt, reqno]{amsart}
\usepackage{amsfonts,amssymb,amscd,amstext,mathrsfs,color}
\usepackage[dvipsnames]{xcolor}
\usepackage{graphicx}
\usepackage{tikz}
\usepackage{hyperref}
\usepackage{mathtools}

\DeclareMathOperator{\codim}{codim}

\numberwithin{equation}{section}

\newtheorem{theorem}{Theorem}[section]
\newtheorem{corollary}[theorem]{Corollary}
\newtheorem{lemma}[theorem]{Lemma}
\newtheorem{proposition}[theorem]{Proposition}

\theoremstyle{definition}
\newtheorem{remark}[theorem]{Remark}

\newtheorem{definition}[theorem]{Definition}

\newcommand{\C}{\mathbb{C}}
\newcommand{\D}{\mathbb{D}}
\newcommand{\B}{\mathbb{B}}

\newcommand{\R}{\mathbb{R}}

\renewcommand{\H}{\mathbb{H}}

\newcommand{\CH}{{\operatorname{CH}}}

\begin{document}
\title{Precise estimates for certain distances in $\R^d$}

\author{Matteo Fiacchi, Nikolai Nikolov} 

\address{M. Fiacchi\\
University of Ljubljana\\
Faculty of Mathematics and Physics\\
Jadranska 19\\
1000 Ljubljana, Slovenia}

\address{Institute of Mathematics, Physics and Mechanics\\
Jadranska 19\\
1000 Ljubljana, Slovenia}

\email{matteo.fiacchi@fmf.uni-lj.si}

\address{N. Nikolov\\
	Institute of Mathematics and Informatics\\
	Bulgarian Academy of Sciences\\
	Acad. G. Bonchev Str., Block 8\\
	1113 Sofia, Bulgaria}

\address{Faculty of Information Sciences\\
	State University of Library Studies
	and Information Technologies\\
	69A, Shipchenski prohod Str.\\
	1574 Sofia, Bulgaria}

\email{nik@math.bas.bg}



\maketitle

\begin{abstract}
We provide sharp estimates for the intrinsic distances of Finsler metrics with precise boundary estimates. These metrics include the Kobayashi-Hilbert metric near strongly convex points, the minimal metric near convex and strongly minimally convex points, and the $k$-quasi hyperbolic metric in $k$-strongly convex domains.

Finally, we prove a characterization result in convex geometry for the $k$-quasi hyperbolic metric.
\end{abstract}

\section{Introduction}
In their famous paper \cite{BaBo}, Balogh and Bonk proved the Gromov hyperbolicity of the Kobayashi distance $k_D$ in strongly pseudoconvex domains, establishing the following estimates for the Kobayashi distance: if $D\subset\C^d$ is a strongly pseudoconvex domain, then there exists $B>0$ such that
\begin{equation}\label{bbintro}
g_D(z,w)-B\leq k_D(z,w)\leq g_D(z,w)+B, \ \ \ \ z,w\in D
\end{equation}
where $g_D$ is a function derived from the Carnot-Carathéodory metric on the boundary $\partial D$.
This estimate, while sufficient to prove Gromov hyperbolicity, does not provide useful information about the distance when the points are close to each other.
The estimates (\refeq{bbintro}) were recently improved in \cite{KNO} in the case of strongly pseudoconvex domains with $\mathcal{C}^{2,\alpha}$-smooth boundary.

A similar situation arises in the real case, where the first named author in \cite{FiaZ} proved the Gromov hyperbolicity of the minimal distance $\rho_D$ (analogous to the Kobayashi metric in the theory of minimal surfaces) in strongly minimally convex domains by showing estimates similar to those of Balogh and Bonk: if $D\subset\R^d$ is strongly minimally convex domain, then there exists $A>1$ and $B>0$ such that
$$A^{-1}d_D(x,y)-B\leq \rho_D(x,y)\leq Ad_D(x,y)+B, \ \ \ \ x,y\in D$$
where 
\begin{equation}\label{d^1}
d_D(x,y)=2\log\left(\frac{\|\pi(x)-\pi(y)\|+h(x)\vee h(y)}{\sqrt{h(x)h(y)}}\right),
\end{equation}
$\pi\colon D\to\partial D$ is such that $\|\pi(x)-x\|=d_{Euc}(x,\partial D)$ and $h(x):=\sqrt{d_{Euc}(x,\partial D)}$.
Also in this case, the estimates were sufficient to prove Gromov hyperbolicity.

The main goal of this paper is to improve the previous estimate in a way that can be applied to various distances naturally defined in real domains.
These distances are all defined through a Finsler metric (see Subsection \ref{Finsler} for a brief introduction), with behavior at the boundary that is asymptotic to $1/2\delta_D$
in the normal component and comparable to $1/\delta_D^{1/2}$
in the tangential component, where $\delta_D(\cdot):=d_{Euc}(\cdot,\partial D)$.

Let us now delve into the details.
Let $D\subset\R^d$ be a domain and $\xi\in\partial D$ be a $\mathcal{C}^2$-smooth boundary point, then there exists a neighborhood $U$ of $\xi$ such that $\pi\colon U\to \partial D$ is a well defined $\mathcal{C}^2$-smooth function. Set for $x,y\in D\cap U$
$$a_D(x,y):=\frac{\|\pi(x)-\pi(y)\|+h(x)\vee h(y)}{\sqrt{h(x)h(y)}}-1,$$
and, for $c>0$, consider the quasi-distance defined as
$$d^c_D(x,y)=2\log(1+ca_D(x,y)).$$
Note that $d^1_D$ is exactly the function that appears in (\refeq{d^1}).

Let $x \in D$ and $v \in \mathbb{R}^d$. If $x$ is sufficiently close to a $\mathcal C^2$-smooth boundary point, we can decompose the vector $v$ into its normal component $v_N$ and its tangential component $v_T$ at the boundary point $\pi(x)$.

The main result of the paper is the following.
\begin{theorem}\label{Main}
Let $D\subset\R^d$ $(d\geq 2)$ be a domain and $\xi\in\partial D$ be a $\mathcal{C}^2$-smooth boundary point. Let $F\colon D\times\R^d\to[0,+\infty)$  be a Finsler metric on $D$ and assume there exist a neighborhood $U$ of $\xi$ 
and $0<c_1<C_1$ such that for all $x\in D\cap U$ and $v\in\R^d$
\begin{equation}\label{Festimatemain}
\max\left\{(1-\omega( \delta_D(x)))\frac{\|v_N\|}{2 \delta_D(x)},c_1\frac{\|v_T\|}{ \delta_D(x)^{1/2}}\right\}\leq F(x,v)\leq (1+\omega( \delta_D(x)))\frac{\|v_N\|}{2 \delta_D(x)}+C_1\frac{\|v_T\|}{ \delta_D(x)^{1/2}}
\end{equation}
where $\omega\colon[0,\varepsilon]\to\R$ is a measurable function with $\int_0^{\varepsilon}\frac{\omega(u)}{u}du<+\infty$. Let $d$ be the intrinsic distance of $F$, then there exist a neighborhood $V\subset\subset U$ of $\xi$ and $0<c_2\leq 1\leq C_2$ such that for all $x,y\in D\cap V$
\begin{equation}\label{destimatemain}
d^{c_2}_D(x,y)\leq d(x,y)\leq d^{C_2}_D(x,y).
\end{equation}
\end{theorem}
The estimates of the main theorem are in the spirit of \cite{KNO}, and they are effective regardless of the relative positions of the two points.

In the second part of the paper, we will prove the estimates (\ref{Festimatemain}) for various Finsler metrics
\begin{enumerate}
	\item Kobayashi-Hilbert metric near strongly convex points (Subsection \ref{subsectKobHilb});
	\item minimal metric near convex and strongly minimally convex points (Subsection \ref{subsectMinimal});
	\item $k$-quasi hyperbolic metric in $k$-strongly convex domains (Subsection \ref{sectkquasihyp}).
\end{enumerate}
As a consequence, we obtain the estimates (\ref{destimatemain}) for the associated intrinsic distances.

The paper concludes with a rigidity theorem for $k$-quasi-hyperbolic metrics (Theorem \ref{theoremappendix}).\\
\newline
\textbf{Acknowledgements.} The first author was partially supported by the European Union (ERC Advanced grant HPDR, 101053085 to Franc Forstneri\v c) and the research program P1-0291 from ARIS, Republic of Slovenia. The second author was partially supported by the Bulgarian National Science Fund, Ministry of Education
and Science of Bulgaria under contract KP-06-N52/3. The first author acknowledges the Institute of Mathematics and Informatics for the hospitality during the visit in November 2024, during which this work was prepared.

\section{Preliminaries}
\textbf{Notations:}\begin{itemize}
	\item For $x\in\R^d$ let $||x||$ denote the standard Euclidean norm of $x$.
	\item For $u,v\in\R^d$ let $\langle u,v\rangle$ denote the Euclidean scalar product of $\R^d$.
	\item Let $\D:=\{x\in\R^2: ||x||<1\}$ be the unit disk in $\R^2$.
	\item For $x,y\in\R^d$ let $$(x,y):=\{tx+(1-t)y: t\in(0,1)\}$$ be the open segment between $x$ and $y$, and $$[x,y]:=\{tx+(1-t)y: t\in[0,1]\}$$ be the closed segment between $x$ and $y$.
	\item For $x\in\R^d$ and $r>0$ we denote with
	$$B(x,r):=\{y\in\R^d:\|x-y\|<r \}$$
	the Euclidean ball with center $x$ and radius $r$.
	\item For $A,B\subseteq\R^d$ nonempty let denote
	$$d_{Euc}(A,B):=\inf\{\|x-y\|: x\in A, y\in B \},$$
	the Euclidean distance between $A$ and $B$. If $A=\{x\}$ is a singleton, we simply write $d_{Euc}(x,\cdot):=d_{Euc}(\{x\},\cdot)$.
	\item If $D\subsetneq\R^d$ is a domain and $x\in\R^d$ let
	$$\delta_D(x):=d_{Euc}(x,\partial D),$$
	be the distance to the boundary.
	\item Let $a,b\in\R$, let denote $a\wedge b:=\min\{a,b\}$ and $a\vee b:=\max\{a,b\}$.
\end{itemize}

\subsection{Finsler metric}\label{Finsler}
In this section, we recall the definition of the Finsler metric.

Let $D\subset\R^d$ be a domain, a \textit{Finsler metric} is a function $F\colon D\times\R^d \simeq TD\to [0,+\infty)$ with the following properties
\begin{enumerate}
	\item $F$ is upper semicontinuous on $D\times\R^d$;
	\item For all $x\in D$, $v\in\R^d$ and $t\in\R$
	$$F(x,tv)=|t|F(x,v).$$
\end{enumerate}

Given a Finsler metric $F$ on $D$ and a piecewise $\mathcal{C}^1$-smooth 
curve $\gamma\colon [a,b]\to D$, we can define the length of $\gamma$ with respect to the metric $F$
$$\ell_F(\gamma):=\int_a^bF(\gamma(t),\dot\gamma(t))dt.$$
Finally, the \textit{intrinsic pseudodistance} associated with $F$ is defined as 
$$d(x,y)=\inf_\gamma\ell_F(\gamma), \ \ \ x,y\in D$$
where the infimum is over all piecewise $\mathcal{C}^1$-smooth curve $\gamma\colon[0,1]\to D$ with $\gamma(0)=x$ and $\gamma(1)=y$.

\subsection{$\mathcal{C}^2$-smooth boundary point}
In this section, we will review some results on the geometry of the domain near a $\mathcal{C}^2$-smooth boundary point.

Let $D\subsetneq\R^d$ be a domain and $\xi\in\partial D$ be a boundary point. We way that $\xi$ is a \textit{$\mathcal{C}^2$-smooth boundary point} if there exists a $\mathcal{C}^2$-smooth local defining function $\rho$ near $\xi$, i.e., there exists a neighborhood $U$ of $\xi$ and $\rho\colon U\to\R$ a $\mathcal{C}^2$-smooth function, such that
\begin{enumerate}
	\item $D\cap U=\{x\in U: \rho(x)<0\}$;
	\item $\nabla \rho\neq0$ in $\partial D\cap U=\{x\in U: \rho(x)=0\}.$
\end{enumerate}
The vector $$n_\xi:=\frac{\nabla\rho(\xi)}{\|\nabla\rho(\xi)\|}$$ is called \textit{unit outer normal} at $\xi$. It does not depend on the choice of the $\mathcal{C}^2$-smooth local defining function $\rho$.

We now state the local version of a well-known lemma.

\begin{lemma}\label{lemma2.1}
	Let $D\subset\R^d$ $(d\geq2)$ be a domain and $\xi_0\in\partial D$ be a $\mathcal C^2$-smooth boundary point. Then there exists $r,\varepsilon>0$ such that if we set
	$$U:=\{x\in\R^d: d_{Euc}(x,\partial D\cap B(\xi_0,r))<\varepsilon\}$$
	we have
	\begin{enumerate}
		\item every point $\xi\in\partial D\cap U$ is a $\mathcal C^2$-smooth boundary point;
		\item for every $x\in U$ there exists a unique $\pi(x)\in\partial D\cap U$ with $||x-\pi(x)||= \delta_D(x)$;
		\item the \textit{signed distance} to the boundary $\rho\colon \R^d\to\R$ given by
		$$\rho(x):=\begin{cases*}
			- \delta_D(x) &\text{if\ $x\in D$}\\
			\delta_D(x) &\text{if\ $x\notin D$}
		\end{cases*}$$is $\mathcal{C}^2$-smooth on $U$;
		\item the fibers of the map $\pi\colon U\to \partial D\cap U$ are
		$$\pi^{-1}(\xi)=\{\xi+tn_\xi: |t|<\varepsilon\}$$
		where $n_\xi$ is the outer unit normal vector of $\partial D$ at $\xi\in\partial D\cap U$;
		\item the gradient of $\rho$ satisfies  for all $x\in U$
		$$\nabla\rho(x)=n_{\pi(x)};$$
		\item the projection map $\pi\colon U\to \partial D$ is $\mathcal{C}^1$-smooth.
		\item if $\gamma\colon [0,1]\to D\cap U$ is a $\mathcal{C}^1$-smooth curve and $\alpha:=\pi\circ\gamma\colon [0,1]\to\partial D\cap U$ its projection to the boundary, then for all $t\in[0,1]$
		$$\frac{1}{2}||(\dot{\gamma}(t))_T||\leq||\dot{\alpha}(t)||\leq 2||(\dot{\gamma}(t))_T||.$$
	\end{enumerate}
\end{lemma}
\proof
See \cite[Lemma 2.1]{BaBo} and \cite[Lemma 4.1]{FiaZ}.
\endproof

Now we have the preliminary results to give a precise meaning to the decomposition into the normal and tangential parts of a vector presented in (\refeq{Festimatemain}).

Let $D\subset\R^d$ ($d\geq2$) be a domain and $\xi\in\partial D$ be a $\mathcal C^2$-smooth boundary point. Let $U$ be the neighborhood of $\xi$ given by the previous Lemma. Then for $x\in D\cap U$ and $v\in\R^d$, we consider the orthogonal decomposition of $v=v_N+v_T$ at (unique) projection point $\pi(x)\in\partial D\cap U$, where
$$v_N:=\langle v,n_{\pi(x)}\rangle n_{\pi(x)}, \ \ \ \  v_T:= v-v_N.$$

\section{Proof of the main theorem}

In this section, we will prove the main result, namely Theorem \ref{Main}. 
In order to make the reading smoother, we will denote all multiplicative constants by $A$, without distinguishing them with other symbols. The same will be done for additive constants, denoted by $B$.

We begin by recalling the function $d_D^c$
mentioned earlier in the introduction.

Let $D\subsetneq\R^d$ be a domain and $\xi\in\partial D$ a $\mathcal{C}^2$-smooth boundary point. Let $U$ be the neighborhood of $\xi$ given by the Lemma \ref{lemma2.1}. For $x,y\in U$, set
$$a_D(x,y):=\frac{\|\pi(x)-\pi(y)\|+h(x)\vee h(y)}{\sqrt{h(x)h(y)}}-1$$
where $h(\cdot):= \delta_D(\cdot)^{1/2}$.
Let $c>0$ and define \begin{equation}\label{dcDdef}d^c_D(x,y)=2\log(1+ca_D(x,y)).\end{equation}
We begin by listing the properties of these functions.
\begin{proposition}\label{propertiesdc}Let $d_D^c$ be the function defined in $(\refeq{dcDdef})$. Then
	\begin{enumerate}
		\item If $0<c_1\leq c_2$, then for all $x,y\in U$
		$$d_D^{c_1}(x,y)\leq d_D^{c_2}(x,y)\leq c_1^{-1}c_2d_D^{c_1}(x,y);$$
		\item If $c\geq1$, $d_D^c$ is a distance on $U$;
		\item For all $c>0$, $d_D^c$ is a quasi-distance on $U$, i.e., there exists $A\geq1$ such that for all $x,y,z\in U$
		$$d_D^c(x,y)\leq A(d_D^c(x,z)+d_D^c(z,y)).$$
	\end{enumerate}
\end{proposition}

In the following lemmas, we will be under the assumptions of Theorem \ref{Main}. Without loss of generality, we can assume that the neighborhood $U$
is the same as in Lemma \ref{lemma2.1}, possibly after shrinking it.

\begin{lemma}\label{lemmavertical} There exists $B>0$ such that for all  piecewise $\mathcal{C}^1$-smooth curve $\gamma\colon [0,1]\to D\cap U$ with endpoints $x,y$ we have
	$$\ell_F(\gamma)\geq\left|\log\left(\frac{h(x)}{h(y)}\right)\right|-B. $$
Moreover, if $\gamma\colon [0,1]\to D\cap U$ is a straight segment contained in $\pi^{-1}(\xi)$ with $\xi\in\partial D\cap U$, then

$$\ell_F(\gamma)\leq\left|\log\left(\frac{h(x)}{h(y)}\right)\right|+B.$$
\end{lemma}
\proof

By (4) in Lemma \ref{lemma2.1} for those $t\in[0,1]$ for which $\dot\gamma(t)$ exists we have $|\frac{d}{dt} \delta_D(\gamma(t))|=||(\dot{\gamma}(t))_N||$, so by the lower estimates of $F$
\begin{align*}\ell_F(\gamma)&\geq \int_0^1(1-\omega( \delta_D(\gamma(t))))\frac{||(\dot{\gamma}(t))_N||}{2 \delta_D(\gamma(t))}dt
	\\&=\int_0^1(1-\omega( \delta_D(\gamma(t))))\frac{|\frac{d}{dt} \delta_D(\gamma(t))|}{2 \delta_D(\gamma(t))}dt
	\\&\geq\left|\int_0^1\frac{\frac{d}{dt} \delta_D(\gamma(t))}{2 \delta_D(\gamma(t))}dt-\frac{1}{2}\int_0^1\frac{\omega( \delta_D(\gamma(t)))\frac{d}{dt} \delta_D(\gamma(t))}{ \delta_D(\gamma(t))}dt\right|
	\\&\geq\left|\log\left(\frac{h(x)}{h(y)}\right)\right|-\frac{1}{2}\int_0^{\varepsilon}\frac{\omega(u)}{u}du
\end{align*}
that implies that statement since by hypothesis $\int_0^{\varepsilon}\frac{\omega(u)}{u}du<+\infty$.

Finally, if $\gamma$ is a straight line segment, using the upper estimates of $F$ we obtain with a similar computation 
$$\ell_F(\gamma)\leq\left|\log\left(\frac{h(x)}{h(y)}\right)\right|+\frac{1}{2}\int_0^{\varepsilon}\frac{\omega(u)}{u}du. $$
\endproof

\begin{lemma}\label{lemmahoriz}
There exists $A>1$ such that for all piecewise $\mathcal{C}^1$-smooth curve $\gamma\colon [0,1]\to D\cap U$ with endpoints $x,y$ we have
$$\ell_F(\gamma)\geq A^{-1}\frac{\|\pi(x)-\pi(y)\|}{ \delta_\gamma^{1/2}}$$
where $\delta_\gamma:=\max_{t\in[0,1]}\delta(\gamma(t))$.
Moreover, for all $x,y\in D\cap U$ with $\delta_D(x)=\delta_D(y)=:\delta_0$ there exists $\gamma\colon [0,1]\to D\cap U$ with endpoints $x,y$ and $ \delta_D(\gamma(t))\equiv\delta_0$ such that 
$$\ell_F(\gamma)\leq A\frac{\|\pi(x)-\pi(y)\|}{\delta_0^{1/2}}.$$
\end{lemma}

\proof
Set $\alpha=\pi\circ\gamma$. By (7) in Lemma \ref{lemma2.1} for those $t\in[0,1]$ for which $\dot\gamma(t)$ exists we have
$$\|\dot\alpha(t)\|\leq2\|(\dot\gamma(t))_T\|.$$
Clearly we have $\int_0^1\|\dot\alpha(t)\|dt\geq\|\pi(x)-\pi(y)\|$, so
\begin{align*}\ell_F(\gamma)&\geq c\int_0^1\frac{\|(\dot\gamma(t))_T\|}{ \delta_D(\gamma(t))^{1/2}}dt
\\&\geq \frac{c}{2}\int_0^1\frac{\|\dot\alpha(t)\|}{ \delta_D(\gamma(t))^{1/2}}dt\\&\geq \frac{c}{2}\frac{\int_0^1\|\dot\alpha(t)\|dt}{\delta_\gamma^{1/2}}\\&\geq\frac{c}{2}\frac{\|\pi(x)-\pi(y)\|}{\delta_\gamma^{1/2}}.\end{align*}

For the second part, since $\partial D\cap U=\partial D\cap B(\xi_0,r)$, the intrinsic and extrinsic distances are bi-Lipschitz, that is there exists $A>1$ such that for all $\xi_1,\xi_2\in D\cap U$ we may find a piecewise $\mathcal C^1$-smooth curve $\alpha\colon[0,1]\to \partial D\cap U$ connecting them with $\int_0^1\|\dot\alpha(t)\|dt\leq A\|\xi_1-\xi_2\|$.
Let $\alpha$ such curve for $\pi(x),\pi(y)\in D\cap U$, and consider $\gamma(t):=\alpha(t)-\delta_0n_{\alpha(t)}$. Clearly $\pi\circ\gamma=\alpha$ and $ \delta_D(\gamma(t))\equiv\delta_0$. Noticing that $(\gamma(t))_N\equiv0$, again by (7) in Lemma \ref{lemma2.1}, we have
\begin{align*}\ell_F(\gamma)&\leq C\int_0^1\frac{\|(\dot\gamma(t))_T\|}{\delta_0^{1/2}}dt
	\\&\leq 2C\frac{\int_0^1\|\dot\alpha(t)\|dt}{\delta_0^{1/2}}\\&\leq 2AC\frac{\|\pi(x)-\pi(y)\|}{\delta_0^{1/2}}.\end{align*}
\endproof

Let us now combine the two previous lemmas.

\begin{lemma}\label{lemmapi}
There exists $B>0$ such that for all piecewise $\mathcal C^1$-smooth curve  $\gamma\colon [0,1]\to D\cap U$ with endpoints $x,y$ we have
$$\ell_F(\gamma)\geq2\log\left(\frac{\|\pi(x)-\pi(y)\|+h(x)\vee h(y)}{\sqrt{h(x)h(y)}}\right)-B.$$
	
\end{lemma}
\proof
The proof is based on a dyadic decomposition similar to the one in \cite[Theorem 1.1]{BaBo}.

If $\pi(x)=\pi(y)$ the estimate follows from Lemma \ref{lemmavertical}, so we may assume $\pi(x)\neq \pi(y)$. 

 Set $H:=\max_{t\in[0,1]} h(\gamma(t))$ and $t_0=\min\{ t\in[0,1]: h(\gamma(t))=H\}$. We divide the curve in two parts, $\gamma_1:=\gamma|_{[0,t_0]}$ and $\gamma_2:=\gamma|_{[t_0,1]}$. 
We will study $\gamma_1$ and $\gamma_2$ separately, starting with the first.

Since $h(x)\leq H$ there exists $k\geq1$ such that
$2^{-k}H<h(x)\leq2^{-(k-1)}H$. Define $0=s_0\leq s_1<\dots s_k\leq t_0$ as follows. Let $s_0=0$ and $s_j=\min\{s\in[0,t_0]:h(\gamma(s))=2^{-(k-j)}H\}$ for $j=1,\dots,k$. Define $x_j=\gamma(s_j)$ for $j=0,\dots,k$.
We divide the problem into two cases.

\textbf{Case 1.1:} there exist $l\in\{1,\dots,k\}$ $$\|\pi(x_{l-1})-\pi(x_l)\|>\frac{2^{-(k-l)}}{8}\|\pi(x)-\pi(y)\|.$$
Since $ \delta_D(\gamma(t))\leq 2^{-(k-l)}H$ for all $t\in [s_{l-1},s_l]$, by Lemma \ref{lemmahoriz}
$$\ell_F(\gamma|_{[s_{l-1},s_l]})\geq A^{-1}\frac{\|\pi(x_{l-1})-\pi(x_l)\|}{2^{-(k-l)}H}\geq A^{-1}\frac{\|\pi(x)-\pi(y)\|}{H}.$$
So by Lemma \ref{lemmavertical}
\begin{align*}\ell_F(\gamma_1)&=\ell_F(\gamma|_{[0,s_{l-1}]})+\ell_F(\gamma|_{[s_{l-1},s_l]})+\ell_F(\gamma|_{[s_{l},t_0]})\\&\geq\log\left(\frac{h(x_{l-1})}{h(x)}\right)+A^{-1}\frac{\|\pi(x)-\pi(y)\|}{H}+\log\left(\frac{H}{h(x_{l})}\right)-B
	\\&\geq\log\left(\frac{H}{h(x)}\right)+A^{-1}\frac{\|\pi(x)-\pi(y)\|}{H}-B.\end{align*}
We retain this estimate and move on to the other case.

\textbf{Case 1.2:} for all $j=1,\dots k$
$$\|\pi(x_{j-1})-\pi(x_j)\|\leq\frac{2^{-(k-j)}}{8}\|\pi(x)-\pi(y)\|.$$
This implies that if we set $t_1=s_k$
$$\|\pi(x)-\pi(\gamma(t_1))\|\leq \sum_{j=1}^{k}\|\pi(x_{j-1})-\pi(x_j)\|\leq \frac{1}{4}\|\pi(x)-\pi(y)\|.$$
Moreover, again by Lemma \ref{lemmavertical}
$$\ell_F(\gamma|_{[0,t_1]})\geq\log\left(\frac{H}{h(x)}\right)-B.$$
Now we address $\gamma_2$: reasoning in the same way, we have two cases.

\textbf{Case 2.1:}
$$\ell_F(\gamma_2)\geq\log\left(\frac{H}{h(y)}\right)+A^{-1}\frac{\|\pi(x)-\pi(y)\|}{H}-B.$$

\textbf{Case 2.2:} there exists $t_2\in [t_0,1]$ such that

$$\|\pi(\gamma(t_2))-\pi(y)\|\leq  \frac{1}{4}\|\pi(x)-\pi(y)\|$$
and
$$\ell_F(\gamma|_{[t_2,1]})\geq\log\left(\frac{H}{h(y)}\right)-B.$$
We now need to consider all four combinations. If Case 1.2 and Case 2.2 occur simultaneously, we have 
$$\|\pi(\gamma(t_1))-\pi(\gamma(t_2))\|\geq \|\pi(x)-\pi(y)\|-\|\pi(x)-\pi(\gamma(t_1))\|-\|\pi(\gamma(t_2))-\pi(y)\|\geq\frac{1}{2}\|\pi(x)-\pi(y)\|,$$
so by Lemma \ref{lemmahoriz}
\begin{align*}\ell_F(\gamma)&=\ell_F(\gamma_{[0,t_1]})+\ell_F(\gamma_{[t_1,t_2]})+\ell_F(\gamma_{[t_2,1]})	\\&\geq\log\left(\frac{H}{h(x)}\right)+A^{-1}\frac{\|\pi(\gamma(t_1))-\pi(\gamma(t_2))\|}{H}+\log\left(\frac{H}{h(y)}\right)-B
	\\&\geq2\log\left(\frac{H}{\sqrt{h(x)h(y)}}\right)+A^{-1}\frac{\|\pi(x)-\pi(y)\|}{H}-B.
	\end{align*}
In the other three combinations, it is easy to see that we still obtain 
$$\ell_F(\gamma)\geq2\log\left(\frac{H}{\sqrt{h(x)h(y)}}\right)+A^{-1}\frac{\|\pi(x)-\pi(y)\|}{H}-B.$$
Now the minimum of the function $$u\mapsto 2\log\left(\frac{u}{\sqrt{h(x)h(y)}}\right)+ A^{-1}\frac{\|\pi(x)-\pi(y\|}{u}$$ is at $u=\frac{A^{-1}}{2}\|\pi(x)-\pi(y)\|>0$, so
$$\ell_F(\gamma)\geq 2\log\left(\frac{\|\pi(x)-\pi(y)\|}{\sqrt{h(x)h(y)}}\right)-B.$$
Finally, combining with Lemma \ref{lemmavertical} we have
\begin{align*}
\ell_F(\gamma)&\geq\max\left\{ 2\log\left(\frac{\|\pi(x)-\pi(y)\|}{\sqrt{h(x)h(y)}}\right),2\log\left(\frac{h(x)\vee h(y)}{\sqrt{h(x)h(y)}}\right)\right\}-B
\\&=	2\log\left(\max \left\{\frac{\|\pi(x)-\pi(y)\|}{\sqrt{h(x)h(y)}}, \frac{h(x)\vee h(y)}{\sqrt{h(x)h(y)}}\right\} \right)-B
\\&\geq 2\log\left(\frac{\|\pi(x)-\pi(y)\|+h(x)\vee h(y)}{\sqrt{h(x)h(y)}}\right)-B.\end{align*}

\endproof

From the previous lemma, we easily obtain the following corollary.

\begin{corollary}\label{corokpoint}
For all $V\subset\subset U$, there exists $B>0$ such that for all $x\in D\cap V$
$$\inf_{y\in D\backslash U}d(x,y)\geq\frac{1}{2}\log\left(\frac{1}{ \delta_D(x)}\right)-B.$$
\end{corollary}	
\proof
Let $\gamma\colon [0,1]\to D$ be a piecewise $\mathcal{C}^1$-smooth curve with $\gamma(0)=x\in D\cap V$ and $\gamma(1)\in D\backslash U$. Consider
$t^*:=\inf\{t\in [0,1]: \gamma(t)\in D\backslash U\}$ and set $y=\gamma(t^*)$, then by Lemma \ref{lemmapi} there exists $B>0$ such that
\begin{align*}\ell_F(\gamma)\geq\ell_F(\gamma|_{[0,t^*]})&\geq 2\log\left(\frac{\|\pi(x)-\pi(y)\|+h(x)\vee h(y)}{\sqrt{h(x)h(y)}}\right)-B
\\&=\frac{1}{2}\log\left(\frac{1}{ \delta_D(x)}\right)+2\log\left(\frac{\|\pi(x)-\pi(y)\|+h(x)\vee h(y)}{\sqrt{h(y)}}\right)-B.\end{align*}
Since $V$ is relatively compact in $U$ we may find $B>0$ such that
$$2\log\left(\frac{\|\pi(x)-\pi(y)\|+h(x)\vee h(y)}{\sqrt{h(y)}}\right)>-B$$
obtaining
$$\ell_F(\gamma)\geq\frac{1}{2}\log\left(\frac{1}{ \delta_D(x)}\right)-B.$$ We conclude taking the infimum over all piecewise $\mathcal C^1$-smooth curves.
\endproof
We are now ready to prove the main theorem.

\proof[Proof of Theorem \ref{Main}]
Let $V$ be a neighborhood of $\xi$ relatively compact in $U$ and let $x,y\in D\cap V$. 
We divide the proof into two cases, depending on whether the points are ''close'' or ''far''.

\textbf{Case 1:} $a_D(x,y)\leq 1$.

We want to prove that there exists $A>1$ such that for all $x\in D\cap U$ and $v\in\R^d$ we have
\begin{equation}\label{disderivdN1}\limsup_{t\to0}\frac{d^1_D(x,x+tv)}{t}\leq AF(x,v).
\end{equation}
Since $\log(1+t)\sim t$ if $t$ is small, we need to study
$$\limsup_{t\to0}\frac{2}{t}\left(\frac{||\pi(x)-\pi(x+tv)||+h(x)\vee h(x+tv)}{\sqrt{h(x)h(x+tv)}}-1\right).$$
Now by (7) in Lemma \ref{lemma2.1},
$$\limsup_{t\to0}\frac{\|\pi(x)-\pi(x+tv)\|}{t}\leq 2\|v_T\|,$$
so
$$\limsup_{t\to0}\frac{2}{t}\left(\frac{||\pi(x)-\pi(x+tv)||}{\sqrt{h(x)h(x+tv)}}\right)\leq 4\frac{\|v_T\|}{ \delta_D(x)^{1/2}}.$$
For the second part, notice that
	$$\frac{h(x)\vee h(x+tv)}{\sqrt{h(x)h(x+tv)}}-1=\frac{\sqrt{h(x)\vee h(x+tv)}-\sqrt{h(x)\wedge h(x+tv)}}{\sqrt{h(x)\wedge h(x+tv)}}.$$
Moreover, by (5) in Lemma \ref{lemma2.1} we have $$\sqrt{h(x)\vee h(x+tv)}-\sqrt{h(x)\wedge h(x+tv)}=|\sqrt{h(x+tv)}-\sqrt{h(x)}|=\frac{1}{4}t\|v_N\|+o(t),$$ and so
$$\lim_{t\to0}\frac{2}{t}\left(\frac{h(x)\vee h(x+tv)}{\sqrt{h(x)h(x+tv)}}-1\right)=\frac{\|v_N\|}{2 \delta_D(x)}.$$
Finally, we can find $A>1$ such that for all $x,y\in D\cap U$
$$\frac{\|v_N\|}{2 \delta_D(x)}+4\frac{\|v_T\|}{ \delta_D(x)^{1/2}}\leq AF(x,v),$$
and so we proved (\refeq{disderivdN1}). By \cite[Theorem 1.3]{Vent}, this implies that for all $x,y\in D\cap U$
$$d_D^1(x,y)\leq Ad(x,y)$$
and so by (1) in Proposition \ref{propertiesdc}
$$d(x,y)\geq A^{-1}d_D^1(x,y)\geq d_D^{A^{-1}}(x,y).$$
For the upper bound, we may assume $h(y)\geq h(x)$. Consider $x'=\pi(x)- h(y)^2n_{\pi(x)}\in D\cap U$. Clearly we have $d(x,x')\leq A\log\left(\frac{h(y)}{h(x)}\right)$, so
$$d(x,y)\leq d(x,x')+d(x',y)\leq A\log\left(\frac{h(y)}{h(x)}\right)+A\frac{\|\pi(x)-\pi(y)\|}{h(y)}.$$
Finally, since for all $t\in[0,1]$ $t\leq\log(1+Ct)$ where $C=e-1$, we have
\begin{align*}d(x,y)&\leq2A\log\left(\sqrt\frac{h(y)}{h(x)}\right)+A\frac{\|\pi(x)-\pi(y)\|}{h(y)}\\&\leq 2A \left(\sqrt\frac{h(y)}{h(x)}-1+\frac{\|\pi(x)-\pi(y)\|}{\sqrt{h(x)h(y)}}\right)\\&=2Aa_D(x,y)\\&\leq2\log(1+ACa_D(x,y))\\&=d^{AC}_D(x,y).\end{align*}

\textbf{Case 2:} $a_D(x,y)>1$.

For the upper bound, let $\varepsilon$ be the constant of Lemma \ref{lemma2.1}. We set
$h=(\|\pi(x)-\pi(y)\|+h(x)\vee h(y))\wedge\varepsilon^{1/2}$ and we consider $x'=\pi(x)-h^2n_{\pi(x)}$, $y'=\pi(y)-h^2n_{\pi(y)}$.
Since $U$ has finite diameter, there exists $A>1$ such that $\|\pi(x)-\pi(y)\|\leq Ah$. Now, by Lemmas \ref{lemmavertical} and \ref{lemmahoriz} we have

\begin{align*}d(x,y)&\leq d(x,x')+d(x',y')+d(y',y)\\&\leq \log\left(\frac{h}{h(x)}\right)+A\frac{\|\pi(x)-\pi(y)\|}{h}+\log\left(\frac{h}{h(y)}\right)\\&\leq2\log\left(\frac{h}{\sqrt{h(x)h(y)}}\right)+A\frac{\|\pi(x)-\pi(y)\|}{h}
	 \\&\leq2\log\left(\frac{\|\pi(x)-\pi(y)\|+h(x)\vee h(y)}{\sqrt{h(x)h(y)}}\right)+B
	 	 \\&=2\log(A+Aa_D(x,y))
	 	 \\&\leq2\log(1+Aa_D(x,y))
	 \\&= d^{A}_D(x,y).
	\end{align*}

For the lower bound, 
let $\gamma\colon[0,1]\to D$ be a $\mathcal C^1$ piecewise curve with endpoints $x,y\in D\cap V$.
Assume again $h(y)\geq h(x)$. Now if $\gamma([0,1])\subset D\cap U$ then by Lemma \ref{lemmapi}
\begin{align*}
\ell_F(\gamma)&\geq2\log (1+a_D(x,y))-B
\\&\geq2\log(1+A^{-1}a_D(x,y))
\\&=d_D^{A^{-1}}(x,y).
\end{align*}

In the other case, set $t_1:=\inf\{t\in[0,1]:\gamma(t)\notin U \}$ and $t_2=\sup\{t\in[0,1]:\gamma(t)\notin U \}$. By definition $\gamma|_{[0,t_1]}\subset D\cap U$ and $\gamma|_{[t_2,1]}\subset D\cap U$. Finally, by Corollary \ref{corokpoint}
\begin{align*}\ell_F(\gamma)&\geq\ell_F(\gamma|_{[0,t_1]})+\ell_F(\gamma|_{[t_2,1]})\\&\geq d(x,\gamma(t_1))+d(\gamma(t_2),y)
\\&\geq 2\log\left(\frac{1}{\sqrt{h(x)h(y)}}\right)-B
\\&\geq 2\log(1+a_D(x,y))-B
\\&\geq2\log(1+A^{-1}a_D(x,y))
\\&= d^{A^{-1}}_D(x,y).\end{align*}
So in both cases, if we take the infimum over all piecewise $\mathcal C^1$-smooth curves connecting $x$ with $y$ we obtain the lower bound
$$d(x,y)\geq d^{A^{-1}}_D(x,y).$$
This completes the proof.
\endproof

We conclude the section by noticing that the estimate (\refeq{destimatemain}) holds in a neighborhood of the boundary if all the boundary points satisfy the assumptions of the main theorem.

\begin{corollary}\label{Maincor}
Let $D\subset\R^d$ be a bounded domain. Let $F\colon D\times\R^d\to[0,+\infty)$  be a Finsler metric on $D$ and assume that for all $\xi\in\partial D$ there exists a neighborhood $U_\xi$ such that $(\refeq{Festimatemain})$ holds. Then there exist a neighborhood $U$ of $\partial D$ and $0<c\leq 1\leq C$ such that for all $x,y\in D\cap U$
$$d^{c}_D(x,y)\leq d(x,y)\leq d^{C}_D(x,y).$$ 
\end{corollary}

\section{Finsler metrics that satisfy the assumptions of Theorem \ref{Main}}

In this section, we prove that several natural metrics in domains of $\R^d$
satisfy the assumptions of Theorem \ref{Main}, including the Kobayashi-Hilbert minimal metric in strongly convex points, the minimal metric in convex and strongly minimally convex points, and the $k$-quasi-hyperbolic metric
metric, recently introduced by Zimmer and Wang in \cite{WanZim}, in 
$k$-strongly convex domains.

The upper bounds for these metrics all arise from the $\mathcal C^2$
regularity of the boundary and the decreasing property with respect to the Beltrami-Klein metric of the ball: let $\B^d\subset\R^d$ be the unit ball, then the \textit{Beltrami-Klein metric} of $\B^d$ is

$$\mathcal{C}\mathcal{K}_{\B^d}(x,v)=\left(\dfrac{(1-||x||^2)||v||^2+|\langle x, v\rangle|^2}{(1-||x||^2)^2}\right)^{1/2}=\left(\dfrac{||v||^2}{1-||x||^2}+\dfrac{|\langle x, v\rangle|^2}{(1-||x||^2)^2}\right)^{1/2}.$$

By composing a translation and a dilation, we obtain the Beltrami-Klein metric for a general Euclidean ball.

\begin{proposition}\label{upperestimate}
Let $D\subset\R^d$ $(d\geq 2)$ be a domain and $\xi\in\partial D$ be a $\mathcal C^2$-smooth boundary point. Let $F\colon D\times\R^d\to [0,+\infty)$ be a Finsler metric on $D$ with the property that if $B$ is an Euclidean ball contained in $D$ then
$$F(x,v)\leq \mathcal{C}\mathcal{K}_B(x,v), \ \ \ \ \forall x\in B\subseteq D,v\in\R^d.$$Then there exists $U$ neighborhood of $\xi$ and $C_1,C_2>0$ such that
$$F(x,v)\leq \left((1+C_1 \delta_D(x))\frac{\|v_N\|^2}{4 \delta_D(x)^2}+C_2\frac{\|v_T\|^2}{ \delta_D(x)}\right)^{1/2}.$$
\end{proposition}

The proposition follows immediately from the following upper bound estimate in the ball.

\begin{lemma}\label{estimateball}
Let $B_r\subset\R^d$ $(d\geq 2)$ be an Euclidean ball of radius $r>0$. Then for all $x\in B_r$ different from the center and $v\in\R^d$ we have
$$\mathcal{C}\mathcal{K}_{B_r}(x,v)\leq\left((1+3r^{-1}\delta_{B_r}(x))\frac{\|v_N\|^2}{4\delta_{B_r}(x)^2}+(1+r^{-1}\delta_{B_r}(x))\frac{1}{2r}\frac{\|v_T\|^2}{\delta_{B_r}(x)}\right)^{1/2}.$$
\end{lemma}

\proof
First of all, $\mathcal{C}\mathcal{K}_{\B^d}$ can be rewritten in the following way if $ x $ is not the origin
$$\mathcal{C}\mathcal{K}_{\B^d}(x,v)=\left(\frac{\|v_N\|^2}{(1-\|x\|^2)^2}+\frac{\|v_T\|^2}{1-\|x\|^2}\right)^{1/2}.$$
So
$$\mathcal{C}\mathcal{K}_{B_r}(x,v)=r^{-1}\mathcal{C}\mathcal{K}_{\B^d}(r^{-1}x,v)=\left(\frac{r^2}{(r+\|x\|)^2}\frac{\|v_N\|^2}{(r-\|x\|)^2}+\frac{1}{r+\|x\|}\frac{\|v_T\|^2}{r-\|x\|}\right)^{1/2}.$$
Now, since $\frac{4r^2}{(r+t)^2}\leq 1+\frac{3}{r}(r-t)$ and $\frac{1}{r+t}\leq \frac{1}{2r}(1+\frac{1}{r}(r-t))$ for $t\in [0,r]$, we have
$$\mathcal{C}\mathcal{K}_{B_r}(x,v)\leq \left((1+3r^{-1}\delta_{B_r}(x))\frac{\|v_N\|^2}{4\delta_{B_r}(x)^2}+(1+r^{-1}\delta_{B_r}(x))\frac{1}{2r}\frac{\|v_T\|^2}{\delta_{B_r}(x)}\right)^{1/2}.$$
\endproof

\proof[Proof of Proposition \ref{upperestimate}]
Let $U$ and $\varepsilon$ be as in Lemma \ref{lemma2.1}. For every $x\in U$, the ball $B(\pi(x)-\varepsilon n_\varepsilon,\varepsilon)$ is internally tangent at $\pi(x)$, and moreover $x\in B$. Since $\delta_B(x)=\delta_D(x)$, the estimate follows from the assumptions and from Lemma \ref{estimateball}.
\endproof

\subsection{Kobayashi-Hilbert metric}\label{subsectKobHilb}

In \cite{Kob}, Kobayashi introduced a projectively invariant metric in the domains of $\R^d$, generalizing the well-known Hilbert metric in the convex domains of $\R^d$.

Let $I:=(-1,1)$ and $D\subset\R^d$ be a domain. The \textit{Kobayashi-Hilbert metric} is the Finsler metric $k_D\colon D\times\R^d\to [0,+\infty)$ defined as
$$k_D(x,v):=\inf\{1/|r|:f\colon I\to D\ \text{projective map}, f(0)=x, f'(0)=rv\}.$$
The \textit{Kobayashi-Hilbert pseudodistance} $K_D$ of $D$ is the intrinsic distance of $k_D$.

The Kobayashi-Hilbert pseudodistance can also be characterized in the following way: it is the largest pseudodistance such that for every projective map $f\colon I\to D$ we have for all $s,t\in I$
$$K_D(f(s),f(t))\leq H_I(s,t):=\frac{1}{2}\left|\log\left(\frac{s+1}{s-1}\cdot\frac{t+1}{t-1}\right)\right|.$$
Notice that from the definition that the Kobayashi-Hilbert metric has the \textit{decreasing property}, i.e., if $D_1 \subseteq D_2$, then $$k_{D_2}(x,v) \leq k_{D_1}(x,v)$$ for all $x \in D_1$ and $v \in \mathbb{R}^d$.

It turns out that Kobayashi-Hilbert metric has an explicit expression: 
set $$Funk(x,v):=\sup\{t>0: x+t^{-1}v\not\in D\}$$ with the notation $\sup\varnothing=0$, then
$$k_D(x,v)=\frac{1}{2}(Funk(x,v)+Funk(x,-v)).$$

As happens in complex geometry with the Kobayashi and Carathéodory metrics, we can dualize the definition of Kobayashi-Hilbert distance and introduce the \textit{Carathéodory-Hilbert distance} in the following way
$$C_D(x,y):=\sup\{H_I(f(x),f(y)):f\colon D\to I\ \text{projective map}\}.$$

From Schwarz's lemma from projective self-maps on $I$, it follows that $C_D \leq K_D$. Moreover, if we denote by $\hat{D}$ the convex hull of $D$, then
$$C_D=C_{\hat D}.$$
Finally, if $D$ is convex, the two distances coincide and are equal to the Hilbert metric.

In order to obtain estimates for the Kobayashi-Hilbert metric near strongly convex points, we need a good lower bound in the case of the ball. Recall that the Hilbert metric in the ball coincides with the Beltrami-Klein metric.

\begin{lemma}\label{estimateballlower}
	Let $B_r\subset\R^d$ $(d\geq 2)$ be an Euclidean ball of radius $r>0$. Then for all $x\in B_r$ that is not the center and $v\in\R^d$ we have
	$$\mathcal{C}\mathcal{K}_{B_r}(x,v)\geq\left(\frac{\|v_N\|^2}{4\delta_D(x)^2}+\frac{1}{2r}\frac{\|v_T\|^2}{\delta_D(x)}\right)^{1/2}.$$
\end{lemma}
\proof
It easily follows from the calculations in the proof of Lemma \ref{estimateball}.
\endproof

We can now prove the lower estimates at strongly convex points, assuming that the domain is convex and bounded.

We recall that, given a domain $D\subseteq\R^d$, a boundary point $\xi\in\partial D$ is said to be \textit{strongly convex} if, there exists a (or equivalently, for all) $\mathcal{C}^2$-smooth local defining function $\rho\colon U\to \R$ at $\xi$ such that $Hess_\xi(\rho)$, the Hessian of $\rho$ at $\xi$, is positive definite on the tangent space $T_\xi\partial D:=\{v\in\R^d: \langle v,n_\xi\rangle=0 \}$.

\begin{proposition}\label{hilblowerconvex}
	Let $D\subset\R^d$ $(d\geq 2)$ be a bounded convex domain and $\xi\in\partial D$ a strongly convex boundary point. Then there exist a neighborhood $U$ of $\xi$ and $c>0$ such that all $x\in D\cap U$ and $v\in\R^d$
	$$k_D(x,v)\geq \left(\frac{\|v_N\|^2}{4 \delta_D(x)^2}+c\frac{\|v_T\|^2}{ \delta_D(x)}\right)^{1/2}.$$
\end{proposition}
\proof
Since $D$ is bounded convex and $\xi$ is strongly convex, it is contained in a Euclidean ball $B_R$ of radius $R$ tangent at $\xi$. The estimate follows immediately from Lemma \ref{estimateballlower} and the decreasing property of the Kobayashi-Hilbert metric.
\endproof

In order to obtain a lower estimate in the case where we do not assume global convexity, we need the following strong localization result. If $D\subset\R^d$ is a convex domain, we say that a boundary point $\xi\in\partial D$ is \textit{strictly convex} if for all $\eta\in\partial D\backslash\{\xi\}$ the segment $(\xi,\eta)$ is contained in $D$.
Moreover, we define 
$$\delta_{D}(x,v):=d_{Euc}(x,(x+\R v)\cap\partial D).$$

\begin{proposition}[Strong localization of Kobayashi-Hilbert metric] Let $D\subset\R^d$ be a domain and $\xi\in\partial D$ be a boundary point. Assume there exists a neighborhood $U$ of $\xi$ such that $D\cap U$ is convex and $\xi$ is a strictly convex for $D\cap U$. Then there exists a neighborhood $V\subset\subset U$ of $\xi$ and $C>0$ such that for all $x\in D\cap V$ and $v\in\R^d$
\begin{equation}\label{localkobhil}
k_{D\cap U}(x,v)\leq (1+C\delta_D(x))k_D(x,v).
\end{equation}
\end{proposition}
\proof
Let $\mathbb{S}^{d-1}$ denote the unit sphere in $\R^d$. Set $D'=D\cap U$. 
Let $x\in D'$ and $v\in\mathbb{S}^{d-1}$, then if the two intersections of $x+\R v$ with $\partial D'$ are both in $\partial D$, then $k_D(x,v)=k_{D'}(x,v)$, and therefore (\refeq{localkobhil}) holds.

For this reason we focus on the case where at least one of the two intersections is not in $\partial D$.

\textbf{Claim 1:}
 There exists $V_1\subset\subset U$ neighborhood of $\xi$ such that for all $x\in D\cap V_1$ and $v\in\mathbb{S}^{d-1}$ we have $\delta_D(x,v)=\delta_{D'}(x,v).$ 
\proof
By contradiction, there exists a sequence of points $x_n$ in $D$ converging to $\xi$ and a sequence of lines $l_n$ passing through $x_n$ such that their intersections with $\partial D'$ (which we denote by $a_n$ and $b_n$) are both in $\partial U$.

Up to subsequences, we can assume that $a_n\to a_\infty$ and $b_n\to b_\infty$. Notice that $a_\infty, b_\infty$ and $\xi$ are distinct. This means that the segment $[a_\infty,b_\infty]$ is contained in $\partial D'$,  but $\xi\in (a_\infty,b_\infty)$, violating the strict convexity of $\xi$.
\endproof

Let $\pi(x)\in\partial D'$ be a closest point (not necessarily unique) and set $n_x:=\frac{\pi(x)-x}{\|\pi(x)-x\|}$ and $v_N:=|\langle v,n_x\rangle|$.

\textbf{Claim 2:}
There exist $V_2\subseteq V_1$ neighborhood of $\xi$ and $c_1>0$ such that if $x\in D\cap V_2$ and $v\in\mathbb{S}^{d-1}$, then if at least one of the two intersections of $x+\R v$ with $\partial D'$ is not in $\partial D$, then $v_N\geq c_1>0$.
\proof
By contradiction, there exists $x_n\to\xi$ and $v_n\in\mathbb{S}^{d-1}$ with $(v_n)_N\to0$ such that one of the intersections of $ l_n:=x_n+\R v_n$ with $\partial D'$ is not in $\partial D$ (let us denote it by $a_n$). Up to subsequence, we can suppose that $a_n\to a_\infty\in\partial D'$ and $v_n\to v_\infty\in \mathbb{S}^{d-1}$. In this way, $ l_n$ converges to a line $ l_\infty$ given by $\xi+\R v_\infty$. Notice that $ l_\infty$ is a line tangent at $\xi$, since is the limit of tangent lines $\pi(x_n)+\R(v_n-(v_n)_N)$. Finally, since $ a_\infty\in l_\infty$ and $ a_\infty\neq\xi$, the segment $[\xi, a_\infty]
$ is contained in $\partial D'$, violating the strict convexity of $\xi$.
\endproof

\textbf{Claim 3:}
For all $x\in V_2$ and $v\in\mathbb{S}^{d-1}$ we have
$$\frac{\delta_{D'}(x)}{\delta_{D'}(x,v)}\geq v_N.$$
\proof
Consider the half-space tangent at $\pi(x)$
$$H:=\{y\in\R^d:\langle y-\pi(x),x-\pi(x) \rangle>0\}.$$ 
Since $D'\subset H$ and $\delta_{D'}(x)=\delta_H(x)$, we have
$$\frac{\delta_{D'}(x)}{\delta_{D'}(x,v)}\geq\frac{\delta_{H}(x)}{\delta_{H}(x,v)}=v_N.$$
\endproof

Let $x\in V_2$ and $v\in\mathbb{S}^{d-1}$.  If least one of the two intersections of $x+\R v$ with $\partial D'$ is not in $\partial D$, by Claim 2 $v_N>c_1>0$, then if we set $c_2:=d_{Euc}(D\cap V_2, D\backslash U)>0$ we have by Claim 3
\begin{align*}k_{D'}(x,v)&\leq\frac{1}{2}\left(\frac{1}{\delta_{D'}(x,v)}+\frac{1}{c_2}\right)
\\&\leq\frac{1}{2}\left(\frac{1}{\delta_{D'}(x,v)}+\frac{1}{c_1c_2}v_N\right)
\\&\leq\frac{1}{2}\left(\frac{1}{\delta_{D'}(x,v)}+\frac{1}{c_1c_2}\frac{\delta_{D'}(x)}{\delta_{D'}(x,v)}\right)
\\&=(1+C\delta_{D'}(x))\frac{1}{2\delta_{D'}(x,v)}
\\&=(1+C\delta_{D'}(x))\frac{1}{2\delta_{D}(x,v)}
\\&\leq(1+C\delta_{D}(x))k_D(x,v).
\end{align*}
\endproof

\begin{remark}
If, under the assumptions of the previous proposition, we add the further assumption that $D$ is hyperbolic (that is, does not contain affine lines), then we have (\refeq{localkobhil}) for every relatively compact $V$ in $U$ (with the constant $C$ that depends on $V$).
\end{remark}

Finally, by combining the localization result with Proposition \ref{hilblowerconvex}, we obtain.

\begin{proposition}
	Let $D\subset\R^d$ $(d\geq 2)$ be a domain and $\xi\in\partial D$ be a strongly convex boundary point. Then there exist a neighborhood $U$ of $\xi$ and $c_1,c_2>0$ such that all $x\in D\cap U$ and $v\in\R^d$
	$$k_D(x,v)\geq \left((1-c_1 \delta_D(x))\frac{\|v_N\|^2}{4 \delta_D(x)^2}+c_2\frac{\|v_T\|^2}{ \delta_D(x)}\right)^{1/2}.$$
\end{proposition}	

\subsection{Minimal metric}\label{subsectMinimal}

In real Euclidean space, Forstneri\v{c} and Kalaj in \cite{FK} defined the \textit{minimal metric}, the analog of the Kobayashi metric in the theory of minimal surfaces.

A map $f\colon\D\to \R^d$ $(d\geq2)$ is said to be \textit{conformal} if for all $\zeta\in\D$ we have
$$||f_x(\zeta)||=||f_y(\zeta)|| \ \mbox{and} \ \langle f_x(\zeta), f_y(\zeta)\rangle=0$$
where $\zeta=(x,y)$ are the coordinates of $\D\subset\R^2$. Moreover, we say that $f$ is \textit{harmonic} if every component of $f$ is harmonic.
If $D\subset\R^d$ $(d\geq3)$ is a domain we denote by $\CH(\D, D)
$ the space of conformal harmonic maps $f\colon\D\to D$.

The \textit{minimal metric} of $D$ is the Finsler metric given by
$$g_D(x,v)=\inf\{1/r: f\in \CH(\D,D), f(0)=x, f_x(0)=rv \}, \ \ \ x\in D, v\in\R^d,$$
and the associated intrinsic distance $\rho_D\colon D\times D\to [0,+\infty)$ is called \textit{minimal pseudodistance}.

As is clearly evident from the definition, the minimal metric also satisfies the decreasing property.

We now introduce the concept of a strongly minimal convex point, which is the analogue in minimal surface theory of strongly convex points in Hilbert geometry and strongly pseudoconvex points in complex analysis.

\begin{definition}[Minimal strongly convex]
Let $D\subset\R^d$ $(d\geq 3)$ be a domain. The boundary point $\xi\in\partial D$ is called \textit{strongly minimally convex} if there exists (or equivalently, for all)  $\mathcal{C}^2$-smooth local defining function $\rho\colon U\to \R$ at $\xi$ such that the smallest two eigenvalues $\lambda_1$ and $\lambda_2$ of $Hess_\xi(\rho)|_{T_{\xi}\partial D}$ satisfies
$$\lambda_1+\lambda_2> 0.$$
\end{definition}

Note that a strongly convex boundary point is strongly minimally convex.

For bounded strongly minimally convex domains (i.e., every boundary point is strongly minimally convex), Drinovec-Drnov\v{s}ek and Forstneri\v{c} have shown in  \cite{DrFor} that we have the following lower estimate for the minimal metric: let $D\subset\R^d$ bounded strongly minimally convex, then there exist $c_1>0$ such that for all $x\in D$ and $v\in\R^d$
\begin{equation}\label{DrFortg}
g_D(x,v)\geq c_1\frac{\|v\|}{ \delta_D(x)^{1/2}}.
\end{equation}
Moreover, there exist $c_2>0$ such that for all $x\in D$ close to $\partial D$ and $v\in\R^d$ we have
$$g_D(x,v)\geq c_2\frac{\|v_N\|}{ \delta_D(x)}.$$

The previous estimates were sufficient for the first author to prove the Gromov hyperbolicity of bounded strongly minimally convex domains (see \cite{FiaZ}). However, for the second estimate (in the normal direction), it is not clear how to improve it to fall within the hypothesis of Theorem \ref{Main}.
For this reason, we restricted ourselves to the convex and locally convex case.

Note that at a \textit{locally convex} boundary point $\xi\in\partial D$, that is, there exists a neighborhood $U$ of $\xi$ such that $D\cap U$ is convex, $\xi$ is strongly minimally convex if and only if 
then $Hess_\xi(\rho)|_{T_\xi\partial D}$ has at least $d-2$ positive eigenvalues, i.e., it is 2-strongly convex in the sense of Definition \ref{kstrcvx}.

\begin{proposition}\label{lowerminiaml}
	Let $D\subset\R^d$ $(d\geq 3)$ be a convex domain and $\xi\in\partial D$ be a point  strongly minimally convex, then there exist a neighborhood $U$ of $\xi$ and $c>0$ such that for all $x\in D\cap U$ and $v\in\R^d$ we have
	$$g_D(x,v)\geq\max\left\{\frac{\|v_N\|}{2 \delta_D(x)},c_2\frac{\|v_T\|}{ \delta_D(x)^{1/2}}\right\}.$$
\end{proposition}	
\proof
The tangential component easy follows from (\refeq{DrFortg}) and the localization result \cite[Theorem 8.5]{DrFor}.

For the normal component, let $x\in D$ close to $\xi$, and consider the half-space
$$H:=\{x\in\R^d: \langle x-\pi(x),n_{\pi(x)}\rangle<0\}.$$
By convexity $D\subset H$, the decreasing property of the minimal metric and \cite[Lemma 5.3]{FiaZ} 
$$g_D(x,v)\geq g_H(x,v)=\frac{\|v_N\|}{2 \delta_D(x)}.$$
\endproof

Finally, from the strong localization result in \cite[Theorem 8.5]{DrFor}, we can slightly relax the convexity assumption, passing from global to local convexity.

\begin{proposition}
	Let $D\subset\R^d$ be a domain and $\xi\in\partial D$ be a strongly minimally convex and locally convex boundary point. Then there exist a neighborhood $U$ of $\xi$ and $c_1,c_2>0$ such that for all $x\in D\cap U$ and $v\in\R^d$ we have
	$$g_D(x,v)\geq\max\left\{(1-c_1 \delta_D(x))\frac{\|v_N\|}{2 \delta_D(x)},c_2\frac{\|v_T\|}{ \delta_D(x)^{1/2}}\right\}.$$
\end{proposition}	

\subsection{$k$-quasi-hyperbolic metric}\label{sectkquasihyp}

In this last subsection, we study a metric recently introduced by Wang and Zimmer in \cite{WanZim}. Let $D\subset\R^d$ be a domain and let $k\in\{1,\dots,d\}$. For all $x\in D$ and $v\in\R^d$, define
$$ \delta_D^{(k)}(x,v):=\sup\{d_{Euc}(x,(x+V)\cap\partial D): V\subseteq\R^d \ \text{a $k$-dimensional linear subspace with}\  v\in V \}.$$
Then the \textit{(generalized) $k$-quasi-hyperbolic metric} on $D$ is defined by
$$q^{(k)}_D(x,v):=\frac{\|v\|}{2 \delta_D^{(k)}(x,v)},$$
where $x\in D$ and $v\in\R^d$. We denote with $d^{(k)}_D$ the associated intrinsic distance.

Note that, compared to the initial definition by Wang and Zimmer, we normalize the metric with a multiplicative factor $\frac{1}{2}$, since our goal is to prove estimates as in the hypothesis of Theorem \ref{Main}.

The $k$-quasi-hyperbolic metric has a strong relationship with several important metrics
\begin{enumerate}
\item $q^{(d)}_D$ is the quasi-hyperbolic metric;
\item $q^{(1)}_D$ is bi-Lipschitz to the Kobayashi-Hilbert metric;
\item $D$ is convex and $d\geq3$, $q^{(2)}_D$ is bi-Lipschitz to the the minimal metric (see \cite[Proposition 10.1]{WanZim}).

\end{enumerate}
Let us now introduce the natural domains where we study the $k$-quasi-hyperbolic metric.

\begin{definition}[$k$-strongly convex domains]\label{kstrcvx}
Let $D\subset\R^d$ $(d\geq 2)$ be a domain. The boundary point $\xi\in\partial D$ is called \textit{$k$-strongly convex boundary point} if it is locally convex and if there exists (or equivalently, for all)  $\mathcal{C}^2$-smooth local defining function $\rho\colon U\to \R$ at $\xi$ such that $Hess_\xi(\rho)|_{T_\xi\partial D}$ has at least $d-k$ positive eigenvalues.
\end{definition}

\begin{remark}
The condition on the Hessian mentioned above is equivalent to requiring that the local defining function is strongly $k$-plurisubharmonic in the sense of \cite{HL,For}.
\end{remark}

Notice that a boundary point is $1$-strongly convex boundary point if and only if is strongly convex, and it is $2$-strongly convex if and only if it is locally convex and strongly minimally convex.

The following property can be viewed as the real and $k$-dimensional analogue of $2$-convexity in the sense of Mercer \cite{Mer}.

\begin{proposition}\label{pkstrcvx}
Let $D\subset\R^d$ $(d\geq 2)$ be a convex domain and $\xi\in\partial D$ be $k$-strongly convex boundary point, then there exist a neighborhood $U$ of $\xi$ and $C>0$ such that for al $x\in D\cap U$, $v\in\R^d$ 
$$ \delta_D^{(k)}(x,v)\leq C \delta_D(x)^{1/2}.$$
\end{proposition}
\proof 
Since $k$-strong convexity is an open condition, we can find a neighborhood $U$ of $\xi$ such that $\pi(x)$ is $k$-strongly convex for all $x\in U$.
Let $x\in D\cap U$. By the $k$-strong convexity, there exist $\lambda_1(x),\dots,\lambda_{d-k}(x)$ positive eigenvalues and $v_1(x),\dots, v_{d-k}(x)$ orthonormal eigenvectors of $Hess_{\pi(x)}(\rho)|_{T_{\pi(x)}\partial D}$. By compactness we can find $\lambda>0$ such that $\lambda<\lambda_j(x)$ for all $x\in D\cap U$ and $j=1,\dots,d-k$. This means that
$$D\subseteq D_{\pi(x)}:=\left\{y\in\R^d:\langle y-\pi(x),n_{\pi(x)}\rangle+\frac{\lambda}{2}\sum_{j=1}^{d-k}\langle y-\pi(x),v_j(x)\rangle^2<0\right\}.$$
Notice that $\delta_D(x)=\delta_{D_{\pi(x)}}(x)$. Finally, a simple calculation shows that there exists $C>0$ such that for all $v\in\R^d$ $$\delta^{(k)}_{D_{\pi(x)}}(x,v)\leq C\delta_{D_{\pi(x)}}(x)^{1/2},$$ and so
$$\delta^{(k)}_D(x,v)\leq\delta^{(k)}_{D_{\pi(x)}}(x,v)\leq C\delta_{D_{\pi(x)}}(x)^{1/2}=C\delta_D(x)^{1/2}.$$
\endproof

Let us now study the case of the Euclidean ball, showing that all the metrics coincide (except for the quasi-hyperbolic one).

\begin{remark}
Let $D\subset\R^d$ $(d\geq 2)$ be a domain, then by definition for all $1\leq k_1\leq k_2\leq d$, we have for $x\in D$ and $v\in\R^d$
$$q_D^{(k_1)}(x,v)\leq q_D^{(k_2)}(x,v).$$
\end{remark}

\begin{proposition}\label{qkball}
Let $d\geq 2$. For all $k\in\{1,\dots,d-1\}$ we have
$$q^{(k)}_{\B^d}=q^{(1)}_{\B^d}.$$
\end{proposition}
\proof
Since for all $x\in\B^d$ and $v\in\R^d$ $$q_{\B^d}^{(1)}(x,v)\leq\cdots\leq q_{\B^d}^{(k)}(x,v)\leq\cdots\leq q_{\B^d}^{(d-1)}(x,v),$$it is sufficient to prove that $$q_{\B^d}^{(d-1)}(x,v)\leq q_{\B^d}^{(1)}(x,v).$$

The result is obvious if $v=0$, so we may assume $v\neq0$.
Let $\xi\in\partial \B^d$ be a boundary point such that $\delta_{\B^d}^{(1)}(x,v)=\|x-\xi\|$. Let $ l$ be the line joining $x$ with $\xi$. 

If $ l$ passes through the origin $O$, let $H$ any hyperplane containing $ l$, then $\B^d\cap H$ is a ball of dimension $d-1$ centered at $O$. Since $O$, $x$ and $\xi$ are collinear, we have
$$\delta_{\B^d}^{(d-1)}(x,v)\geq d_{Euc}(x,(x+H)\cap\partial \B^d)=\|x-\xi\|=\delta_{\B^d}^{(1)}(x,v),$$
and so $q_{\B^d}^{(d-1)}(x,v)\leq q_{\B^d}^{(1)}(x,v).$

If, on the other hand, $ l$ does not pass through the origin, let $O'$ the projection of $O$ onto $ l$. Now let $H$ be the affine hyperplane containing $ l$ and orthogonal to $\overrightarrow{OO'}$. Notice that $\B^d\cap H$ is a ball of dimension $d-1$, centered at $O'$. Since $O'$, $x$ and $\xi$ are collinear, we obtain again $q_{\B^d}^{(d-1)}(x,v)\leq q_{\B^d}^{(1)}(x,v).$
\endproof

From an explicit calculation, we can obtain the expression of the $q_{\H^d}^{(k)}$ metrics in the half-space $\H^d$.

\begin{proposition}\label{halfspaceqk}
Let $\H^d:=\{x\in\R^d:x_1>0\}$ $(d\geq 2)$ be the half-space. Then for all $k\in\{1,\dots,d-1\}$, $x=(x_1,\dots,x_d)\in \H^d$ and $v=(v_1,\dots,v_d)\in\R^d$, we have
$$q^{(k)}_{\H^d}(x,v)=\frac{|v_1|}{2x_1}.$$
\end{proposition}

Finally, we can prove the necessary estimates for Theorem \ref{Main} at $k$-strongly convex points.

\begin{proposition}
	Let $D\subset\R^d$ $(d\geq 2)$ be a convex domain and let $\xi\in\partial D$ be a $k$-strongly convex boundary point, then there exist $c,C_1,C_2>0$ and a neighborhood $U$ of $\xi$ such that for all $x\in D\cap U$ and $v\in\R^d$ we have
	$$\max\left\{\frac{\|v_N\|}{2 \delta_D(x)},c\frac{\|v_T\|}{ \delta_D(x)^{1/2}}\right\}\leq q_D^{(k)}(x,v)\leq \left((1+C_1 \delta_D(x)))\frac{\|v_N\|^2}{4 \delta_D(x)^2}+C_2\frac{\|v_T\|^2}{ \delta_D(x)}\right)^{1/2}.$$
\end{proposition}
\proof
\textbf{Upper bound:} Let $B$ be an Euclidean ball. Then by Proposition \ref{qkball} and \cite[Propotision 10.1]{WanZim} for all $x\in B$ and $v\in\R^d$ we have (being careful with the different normalization of $q^{(k)}_D$)
$$q^{(k)}_{B}(x,v)=q^{(2)}_{B}(x,v)\leq \mathcal{C}\mathcal{K}_{B}(x,v).$$ 
So the upper estimate follows from Proposition \ref{upperestimate}.

\textbf{Lower bound:} First of all, by Proposition \ref{pkstrcvx} there exists $c>0$ such that for all $x\in D$ and $v\in\R^d$
$$q_D^{(k)}(x,v)=\frac{\|v\|}{2 \delta_D^{(k)}(x,v)}\geq c\frac{\|v\|}{ \delta_D(x)^{1/2}}.$$
For the normal component, it is sufficient to reason as in Proposition  \ref{lowerminiaml}.
\endproof

The estimates from the main theorem immediately imply the Gromov hyperbolicity of the $k$-quasi hyperbolic metric in \textit{$k$-strongly convex domains}, which are bounded convex domains where all the boundary points are $k$-strongly convex.
For more details on Gromov hyperbolicity, see \cite{BH}.

\begin{corollary}
Let $d\geq 2$ and $k\in\{1,\cdots,d-1\}$. If $D\subset\R^d$ is a $k$-strongly convex domain, then $(D,d^{(k)}_D)$ is Gromov hyperbolic.
\end{corollary}
\proof
From Corollary \ref{Maincor}, there exists a neighborhood $U$ of the boundary such that (\refeq{destimatemain}) holds. This implies that there exists $B>0$ such that for all $x,y\in D\cap U$
$$2\log\left(\frac{\|\pi(x)-\pi(y)\|+h(x)\vee h(y)}{\sqrt{h(x)h(y)}}\right)-B\leq d_D^{(k)}(x,y)\leq 2\log\left(\frac{\|\pi(x)-\pi(y)\|+h(x)\vee h(y)}{\sqrt{h(x)h(y)}}\right)+B.$$
The proof concludes as in \cite[Theorem 1.4]{BaBo} and \cite[Proposition 4.4]{FiaZ}.
\endproof

Note that the Gromov hyperbolicity for the $k$-quasi hyperbolic metric in convex domains has been characterized by Wang and Zimmer in \cite[Theorem 1.5]{WanZim}.

We conclude with a remark.

\begin{remark}
Let $D\subset\R^d$ be a $k_1$-strongly-convex domain, then for all $k_1<k_2\leq d-1$ there exists $A>1$ such that for all $x\in D$ and $v\in\R^d$
	$$q_D^{(k_1)}(x,v)\leq q_D^{(k_2)}(x,v)\leq A q_D^{(k_1)}(x,v).$$
	This is a consequence of the Theorem \ref{Main} and Proposition \ref{propertiesdc}.
\end{remark}

\section{A rigidity result in convex geometry}

In this section, we will characterize convex domains $D\subsetneq\R^d$ where
\begin{equation}\label{eqappendix}
	\delta_D^{(1)}(x,v)=\delta_D^{(d-1)}(x,v), \ \ \ \ \forall x\in D, v\in\R^d.
\end{equation}
Let us begin with some basic notions of convex geometry.

Let $D\subsetneq\R^d$ be a convex domain. Let $a\in\partial D$. An affine hyperplane $H$ passing at $a$ is a \textit{supporting hyperplane} at $a$ if $D\cap H=\varnothing.$
A \textit{normal line} at $a$ is a line passing through $a$, orthogonal to a supporting hyperplane at $a$.
A \textit{face} is the convex subset $\partial D\cap H$, where $H$ is a supporting hyperplane.

Let $a\in\R^d$ and $v\in\R^d$ non zero. A \textit{half-space} is a domain of the form
$$\{x\in\R^d: \langle x-a,v\rangle<0\}.$$

A \textit{slab} is a domain of the form
$$\{x\in\R^d: |\langle x-a,v\rangle|<1 \}.$$

We can now state the main result of this section

\begin{theorem}\label{theoremappendix}
Let $D\subsetneq\R^d$ $(d\geq3)$ be a convex domain. Then $(\refeq{eqappendix})$ holds if and only if $D$ is either
\begin{enumerate}
	\item an Euclidean ball;
	\item a half-space;
	\item a slab.
\end{enumerate}
\end{theorem}
\proof
We already proved that in the Euclidean balls (\refeq{eqappendix}) holds (see Proposition \ref{qkball}). It is not difficult to see that the same is true for the half-spaces and slabs.

For the other direction, let $D\subsetneq \R^d$ be a convex domain such that (\refeq{eqappendix}) holds. We divide the proof in several parts.

First of all notice that if $D$ has at most 2 faces, then it is either a half-space or slab, so we may assume that $D$ has at least 3 faces.

\textbf{Part 1:} Let $a,b\in\partial D$ not on the same face, and let $ l_a$, $ l_b$ be two normal lines at $a$ and $b$ respectively. Then $ l_a$ and $ l_b$ are coplanar.

Since $a$ and $b$ are not in the same face, $(a,b)\subset D$. Let $H_a$ and $H_b$ be two supporting hyperplanes at $a$ and $b$ orthogonal to the lines $ l_a$ and $ l_b$ respectively. Consider $m=\frac{a+b}{2}$ the midpoint between $a$ and $b$, and $v=\overrightarrow{ab}$. By (\refeq{eqappendix}) we have
$$\delta_D^{(1)}(m,v)=\delta_D^{(d-1)}(m,v).$$
which means that there exists an affine hyperplane $H$ passing through $a$ and $b$ such that
$$B\left(m,\frac{|\overrightarrow{ab}|}{2}\right)\cap H\subset D.$$
This implies that $H_a\cap H$ and $H_b\cap H$ are parallel. Moreover, $\overrightarrow{ab}\perp(H_a\cap H)$ and $\overrightarrow{ab}\perp(H_b\cap H)$. Consequently,
$ l_a\perp(H_a\cap H)$ and $ l_b\perp(H_b\cap H)$, so since $\codim(H_a\cap H)=\codim(H_b\cap H)=2$, the two lines $ l_a$ and $ l_b$ are coplanar.

\textbf{Part 2:} Each point on the boundary is an extreme point, meaning that all the faces are singletons.

By contradiction, let $a,b\in\partial D$ be two distinct boundary points in the same face $F$ of $D$. By definition, there exists an hyperplane $H$ such that $F=\partial D\cap H$. Let $l_a$ and $l_b$ be the two lines perpendicular to $H$, passing through $a$ and $b$, respectively. Clearly, $l_a$ and $l_b$ are coplanar, so there exists a plane $\Pi$ that contains them. Now since $d\geq3$ we can find a boundary point $c\in\partial D\backslash (F\cup \Pi)$ and a normal line $ l_c$ at $c$ that is not parallel at $\Pi$. Such a line exists because otherwise $D$ is a cylinder over $\Pi$, and it is not difficult to see that the only cylinders where (\refeq{eqappendix}) holds are half-spaces and slabs.

Since $l_a$ and $l_c$ are coplanar but not parallel, they intersect at one point. The same is true for $l_b$ and $l_c$, and this implies that $l_c$ lies in the plane $\Pi$, leading to a contradiction.

\textbf{Part 3:} All the normal lines intersect at one point $O$. 

Let $a,b\in\partial D$ be two distinct boundary points with normal lines $l_a$ and $l_b$. From Part 1, there exists a plane $\Pi$ that contains $l_a$ and $l_b$. Consider $c\in\partial D\backslash \Pi$ and a normal line $l_c$ at $c$ that is not parallel to $\Pi$. We conclude as in Part 2.

\textbf{Part 4:} Each point $a\in\partial D\backslash\{O\}$ has a unique supporting hyperplane, which implies that $\partial D\backslash\{O\}$ is $\mathcal{C}^1$-smooth hypersurface.

Suppose that there exist two distinct normal lines at $a_0\in\partial D$. This implies $a_0$ is the point $O$ of Part 2, and thus, every point $a\in\partial D\backslash\{O\}$ has a unique supporting hyperplane. A classic result from convex geometry (see, for example, \cite{Busemann}) ensures that $\partial D$ is $\mathcal{C}^1$-smooth.

\textbf{Part 5:} $D$ is an Euclidean ball.

Up to a translation, we can suppose that the point $O$ of Part 2 is the origin. Let $f\colon\R^d\to \R$ the function given by $f(x):=\|x\|$. Let $a\in\partial D\backslash\{O\}$, since the (unique) normal line at $a$ passes through the origin, $$T_a\partial D=\{v\in\R^d:\langle x,a\rangle=0 \}=(\R a)^{\perp}.$$
Now $df_a|_{T_a\partial D}=0$, which means that $f$ is constant in $\partial D\backslash\{O\}$, and so $\partial D$ is a sphere centered at the origin.
\endproof

\end{document}